\RequirePackage{fix-cm}
\documentclass[twocolumn,citesort]{svjour3}          
\smartqed  
\usepackage{graphicx}
\usepackage{amsmath}
\usepackage{algorithm}
\usepackage[noend]{algpseudocode}
\usepackage{bbold}
\usepackage[stable]{footmisc}
\usepackage[rightcaption]{sidecap}
\makeatletter
\def\BState{\State\hskip-\ALG@thistlm}
\makeatother
\newcommand{\pushcode}[1][1]{\hskip\dimexpr#1\algorithmicindent\relax}
\begin{document}

\title{Numerical Algorithm for P\'olya Enumeration Theorem
}


\author{Conrad W. Rosenbrock         \and
        Wiley S. Morgan \and
        Gus L. W. Hart \and
        Stefano Curtarolo \and
        Rodney W. Forcade 
}


\institute{Conrad W. Rosenbrock \and Wiley S. Morgan \and Gus L. W. Hart
\at
              Department of Physics and Astronomy, Brigham Young
              University, Provo, UT, 84602 \\
              Tel.: +1-801-422-7444\\
              \email{gus.hart@gmail.com}
\and
Rodney W. Forcade
\at
           Department of Mathematics, Brigham Young
              University, Provo, UT, 84602, USA
\and
Stefano Curtarolo
\at
           Materials Science, Electrical Engineering,
Physics and Chemistry,
           Duke University, Durham, NC 27708, USA
}

\date{Received: date / Accepted: date}

\maketitle

\begin{abstract}
Although the P\'olya enumeration theorem has been used extensively for
decades, an optimized, purely numerical algorithm for calculating its
coefficients is not readily available. We present such an algorithm
for finding the number of unique colorings of a finite set under the
action of a finite group.

\end{abstract}

\hyphenation{co-lo-rings}

\section{Introduction}
\label{intro}
A common problem in many fields involves enumerating the possible
colorings of a finite set. Applying a symmetry or permutation group
reduces the size of the enumerated set by including only those elements that
are unique under the group action. The P\'olya enumeration theorem
counts the number of unique colorings that should be recovered
\cite{Polya:1987}.
The P\'olya theorem has shown its wide range of applications in a variety of contexts, such as
confirming enumerations of molecules in bioinformatics and chemoinformatics
\cite{Deng:2014}; unlabeled, uniform hypergraphs in discrete mathematics
\cite{Qian:2014}; and photosensitisers in photosynthesis research
\cite{Taniguchi:2014}.

Typical implementations of the counting theorem use Computer
Algebra Systems to symbolically solve the polynomial coefficient
problem. However, despite the
widespread use of the theorem, a low-level numerical implementation
for recovering the number of unique colorings is not readily
available. Although a brute-force calculation of the expansion
coefficients for the P\'olya polynomial is straight-forward to
implement, it is prohibitively slow. For instance, we recently used
such a brute force method to confirm enumeration results for a lattice
coloring problem in solid state physics \cite{Hart:2008}. After profiling performance
on more than 20 representative systems, we found that the brute force
calculation of the P\'olya coefficient took as long as the enumeration
problem itself. Here we demonstrate that the performance can be
improved drastically by exploiting the properties of polynomials. The
improved performance also enables harder P\'olya theorem problems to
be easily solved that would otherwise be 
computationally prohibitive
\footnote{For example, in one test we performed, Mathematica
  required close to 5 hours to compute the coefficient, while our
  algorithm found the same answer in 0.2 seconds.}.

We first briefly describe the P\'olya enumeration theorem in Section
\ref{sec:polya}, followed by the algorithm for calculating the
polynomial coefficients in Section \ref{sec:algorithm}. In the final
Section, we investigate the scaling and performance of the algorithm
both heuristically and via numerical experimentation.

\section{P\'olya Enumeration Theorem}
\label{sec:polya}
Because of extensive literature coverage, we do not derive the P\'olya's
theorem here\footnote{The interested reader may refers to
  Refs.~\cite{Polya:1937,Polya:1987}}. Rather, we just state its main
claims by using a simple example.  

\begin{SCfigure*}
\centering{\includegraphics[width=.7\textwidth]{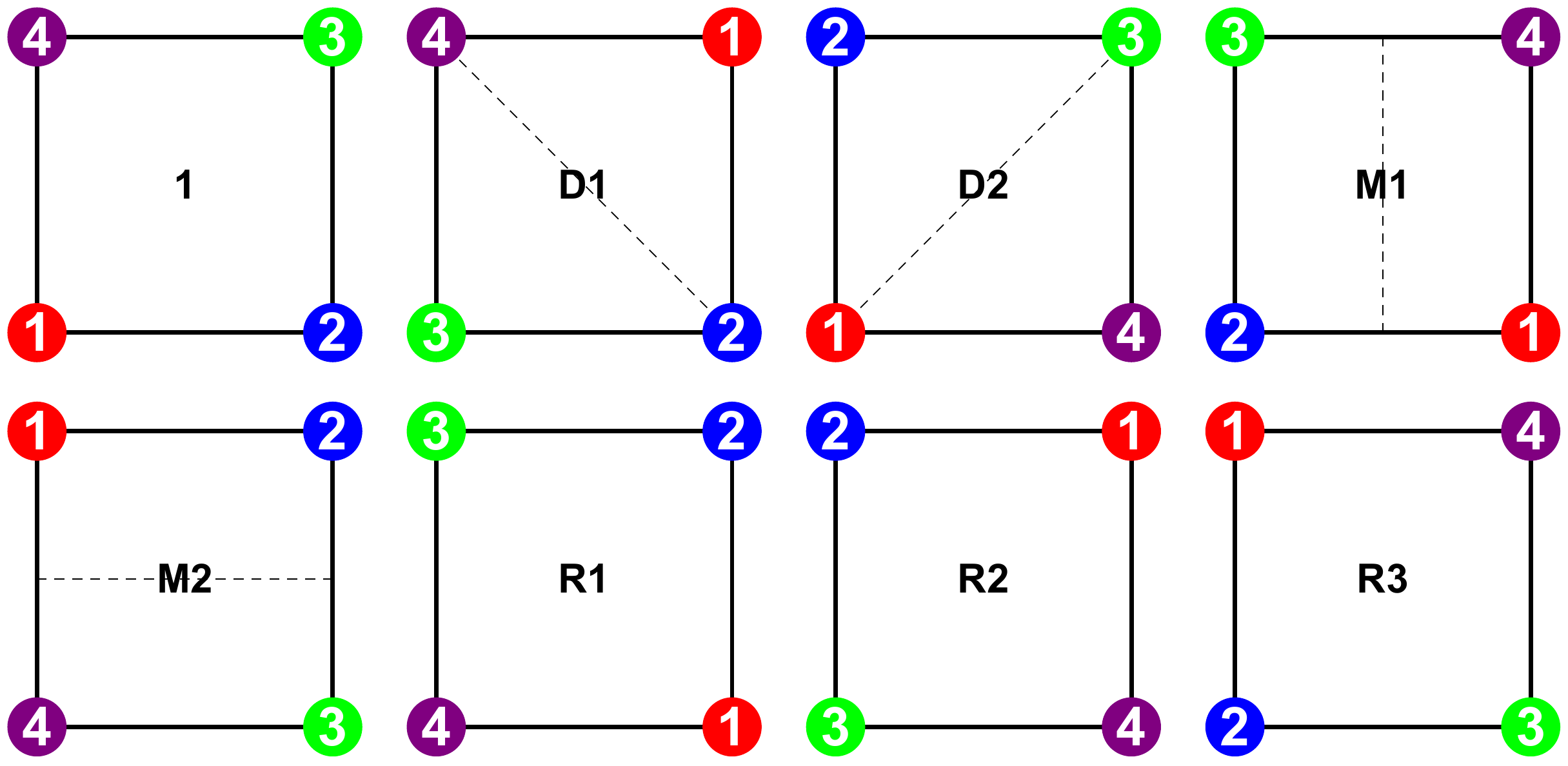}}
\caption[r]{\label{fig:D4GroupOperations} The symmetry group operations
of the square. This group is known as the dihedral group of degree 4,
or D$_4$. The dashed lines are guides to the eye for the horizontal,
vertical and diagonal reflections ({\bf M1},{\bf M2} and {\bf D1},
{\bf D2}).}
\end{SCfigure*}

The square has the set of symmetries displayed in Figure
\ref{fig:D4GroupOperations}. These symmetries include three rotations
(by 90, 180 and 270 degrees; labelled {\bf R1}, {\bf R2}, and {\bf R3}) and four
reflections (one horizontal, one vertical and two for the diagonals;
labelled {\bf M1}, {\bf M2} and {\bf D1}, {\bf D2}). This group is commonly known as the
dihedral group of degree four, or D$_4$ for short\footnote{The
  dihedral groups have multiple, equivalent names. D$_4$ is also
  called Dih$_4$ or the dihedral group of \emph{order} 8 (D$_8$).}.

The group operations of the D$_4$ group can be written in disjoint-cyclic
form as in Table \ref{tab:D4GroupPolynomials}. For each
$r$-cycle in the group, we can write a polynomial in
variables $x_i^r$ for $i=1\dots \xi$, where $\xi$ is the number of colors
used. For this example, we will consider the situation where we want to
color the four corners of the square with just two colors. In that
case we end up with just two variables $x_1, x_2$, which are represented
as $x,y$ in the Table.
 
\begin{table*}[!tb]
\centering
\begin{tabular}{| c | c | c | c | c |}
   \hline
   Op.  & Disjoint-Cyclic & Polynomial & Expanded & Coeff. \\ \hline
$\mathbb{1}$ & $(1)(2)(3)(4)$ & $(x+y)^4$ & $x^4+4 x^3 y+6 x^2 y^2+4 x
y^3+y^4$ & 6 \\
{\bf D1} & $(1,3)(2)(4)$ & $(x^2+y^2)(x+y)^2$ & $x^4+2 x^3 y+2 x^2 y^2+2 x
y^3+y^4$ & 2 \\
{\bf D2} & $(1,2)(3)(4)$ & $(x^2+y^2)(x+y)^2$ & $x^4 + 2 x^3 y + 2 x^2 y^2
+ 2 x y^3 + y^4$ & 2\\
{\bf M1} & $(1,2)(3,4)$ & $(x^2+y^2)^2$ & $x^4 + 2 x^2 y^2 + y^4$ & 2 \\
{\bf M2} & $(1,4)(2,3)$ & $(x^2+y^2)^2$ & $x^4 + 2 x^2 y^2 + y^4$ & 2 \\
{\bf R1} & $(1,4,3,2)$ & $(x^4+y^4)$ & $x^4 + y^4$ & 0\\
{\bf R2} & $(1,3)(2,4)$ & $(x^2+y^2)^2$ & $x^4 + 2 x^2 y^2 + y^4$ & 2\\
{\bf R3} & $(1,2,3,4)$ & $(x^4+y^4)$ & $x^4 + y^4$ & 0\\
   \hline
\end{tabular}
\caption{Disjoint-cyclic form for each group operation in D$_4$ and the
  corresponding polynomials, expanded polynomials and the coefficient
  of the $x^2y^2$ term for each.}
\label{tab:D4GroupPolynomials}
\end{table*}

The P\'olya representation for a single group operation in
disjoint-cyclic form results in a product of
polynomials that we can expand. For example, the group operation {\bf
  D1} has disjoint-cyclic form $(1,3)(2)(4)$ that can be represented
by the polynomial $(x^2+y^2)(x+y)(x+y)$ where the exponent on each
variable corresponds to the length of the $r$-cycle that it
is part of. For a general $r$-cycle, the polynomial takes the form
\begin{equation}
\label{eq:GenericPolynomialOneRCycle}
(x_1^r+x_2^r+\cdots+x_\xi^r),
\end{equation}
\noindent for an enumeration with $\xi$ colors. Most group operations
will have a product of these polynomials for each $r$-cycle in the
disjoint-cyclic form. Once the product of polynomials has been
generated with the group operation, we can simplify it by adding
exponents to identical polynomials. In the example above, $(x+y)(x+y)$
would become $(x+y)^2$; in summary, we exchange the group operations
acting on the set for polynomial representations that obey the
familiar rules for polynomials.

\begin{figure}
\centerline{\includegraphics[width=.35\textwidth]{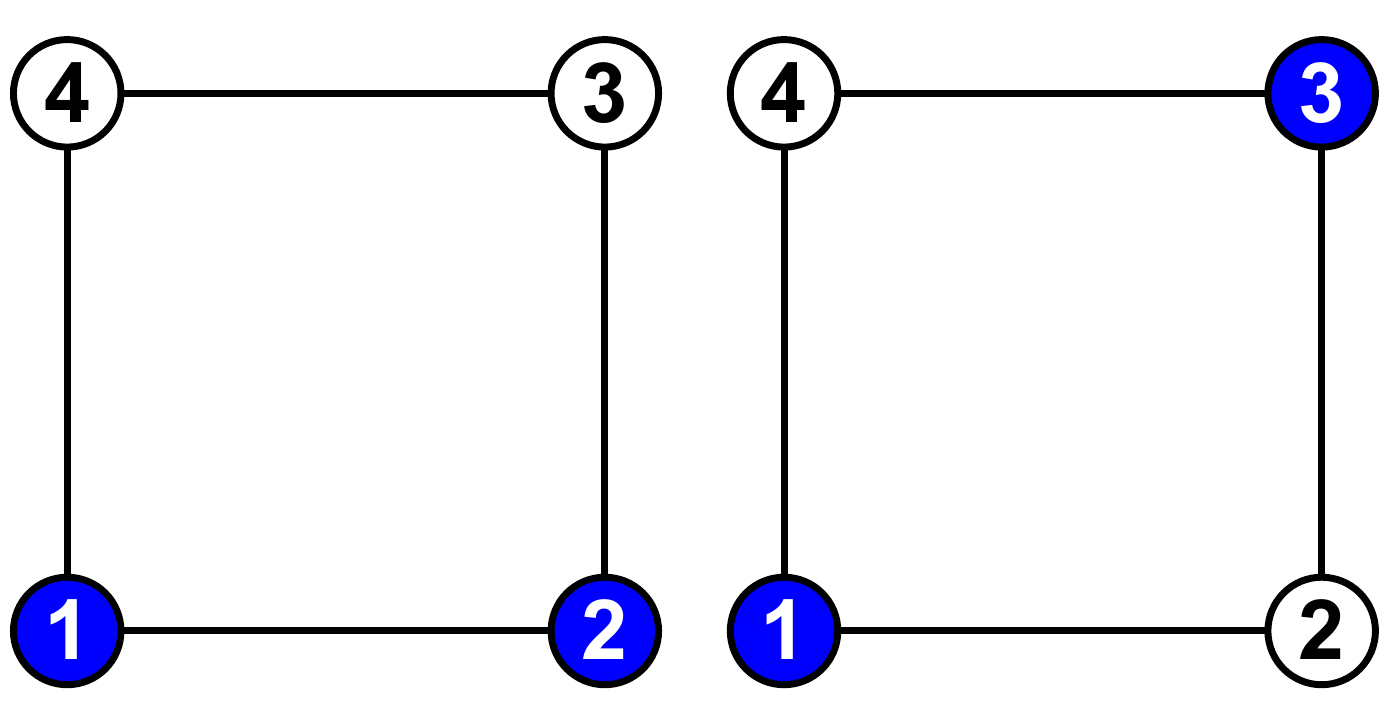}}
\caption[]{\label{fig:UniqueColoringsD4} The two unique ways to color
  the square with two colors and two corners of each color.}
\end{figure}

We will now pursue our example of the possible colorings on the four
corners of the square involving two of each color. Excluding
the symmetry operations, we could come up with ${4 \choose 2} = 6$
possibilities, but some of these are equivalent by symmetry. The
P\'olya theorem will count how many \emph{unique} colorings we should recover. To
find out the expected number of unique colorings, we look at the
coefficient of the term corresponding to the overall color selection
(in this example, two of each color); thus we look for coefficients of
the $x^2y^2$ term for each group operation. These coefficient values
are listed in Table \ref{tab:D4GroupPolynomials}. The sum of these
coefficients, divided by the number of operations in the group, gives
the total number of unique colorings under the entire group action, in
this case $(6+2+2+2+2+0+2+0)/8=16/8=2$. The unique colorings are plotted 
in Figure \ref{fig:UniqueColoringsD4}.

Generally, for a finite set with $F$ elements, and
fixed color concentrations $c_i$ such that $\sum_{i=1}^{\xi} c_i =
F$, the number of unique colorings of the set under the group action
corresponds to the coefficient of the term
\begin{equation}
\label{eq:FixedColorTerm}
T=x_1^{c_1}x_2^{c_2}\dots x_{\xi}^{c_\xi}=\prod_i^{\xi}x_i^{c_i}
\end{equation}
\noindent in the expanded polynomial for each group operation, summed over all
elements in the group. Counting the number of unique colorings at 
fixed concentration amounts to finding the coefficient of a specific
term, known \emph{a priori}, from a product of polynomials.

\section{Coefficient-Finding Algorithm}
\label{sec:algorithm}
We begin by reviewing some well-known properties of polynomials with
respect to their variables. First, for a generic polynomial
\begin{equation}
\label{eq:GenericMultinomial}
(x_1^r+x_2^r+\cdots+x_\xi^r)^d,
\end{equation}
\noindent the exponents of each $x_i$ in the \emph{expanded}
polynomial are constrained to the set
\begin{equation}
\label{eq:PossibleVariableValues}
V=\{0,r,2r,3r,\dots,dr\}.
\end{equation}

Next, we consider the terms in the expansion of the polynomial:
\begin{equation}
\label{eq:MultinomialExpansion}
(x_1^r+x_2^r+\cdots+x_\xi^r)^d = \sum_{k_1,k_2,\dots,k_\xi} \mu_k
\prod_{i=1}^\xi x_i^{rk_i}
\end{equation}

\noindent where the sum is over all possibles sequences $k_1,k_2,\dots,k_\xi$ such
that the sum of the exponents (represented by the sequence in $k_i$)
is equal to $d$,
\begin{equation}
\label{eq:SumOfExponentsCondition}
k_1+k_2+\cdots+k_\xi=d.
\end{equation}

The coefficients $\mu_k$ in the polynomial expansion Equation
(\ref{eq:MultinomialExpansion}) are found using the multinomial
tcoefficients

\begin{eqnarray}
\mu_k & = & {n \choose {k_1,k_2,\dots,k_\xi}} =
\frac{n!}{k_1!k_2!\cdots k_\xi!} \nonumber \\
& = & {{k_1}\choose{k_1}}{{k_1+k_2}\choose{k_2}} \cdots {{k_1+k_2+\cdots+k_\xi}
\choose{k_\xi}} \nonumber \\
\label{eq:MultinomialCoefficientEquation}
& = & \prod_{i=1}^\xi {{\sum_{j=1}^i k_j} \choose k_i}.
\end{eqnarray}

Finally, we define the polynomial
(\ref{eq:GenericPolynomialOneRCycle}) for an arbitrary
group operation $G^i \in \bf{G}$ as\footnote{We will use
  Greek subscripts to label the polynomials in the product and Latin
  subscripts to label the variables within any of the polynomials.}
\begin{equation}
\label{eq:ProductOfMultinomials}
P^i(x_1,x_2,\dots,x_\xi) = \prod_{\alpha=1}^m M_\alpha^{r_{\alpha}}(x_1,x_2,\dots,x_\xi)
\end{equation}

\noindent where each $M_\alpha^{r_{\alpha}}$ is a polynomial for the
$\alpha^{\mathrm{th}}$ distinct $r$-cycle
of the form (\ref{eq:GenericMultinomial}) and $d_\alpha$ is
substituted for the value of
$d$ (which is the multiplicity of that $r$-cycle); $m$ is the
number of distinct values of $r_{\alpha}$ in $P^i$. 

\emph{Since we know the fixed concentration term
$T=\prod_{i=1}^\xi T_i=\prod_{i=1}^\xi x_i^{c_i}$ in advance (see
equation (\ref{eq:FixedColorTerm})), we can limit the
possible sequences of $k_i$ for which multinomial coefficients are
calculated}. This is the key idea of the algorithm and the reason for
its high performance.

For each group operation $G^i$, we have a product of
polynomials $M_\alpha^{r_{\alpha}}$. We begin filtering the sequences
by choosing only those combinations of values $v_{i\alpha} \in
V_\alpha = \{v_{i\alpha}\}_{i=1}^{d_{\alpha}+1}$ for which the sum 
\begin{equation}
\label{eq:SumOfExponentsEqualsTarget}
\sum_{\alpha=1}^m v_{i\alpha} = T_i
\end{equation}

\noindent where $V_\alpha$ is the set from
eqn. (\ref{eq:PossibleVariableValues}) for multinomial $M_\alpha^{r_{\alpha}}$.

We first apply constraint (\ref{eq:SumOfExponentsEqualsTarget}) to the $x_1$
term across the product of polynomials to find a set of values
$\{k_{1\alpha}\}_{\alpha=1}^m$ that could give exponent $T_1$ once all the
polynomials' terms have been expanded. Once a value $k_{1\alpha}$
has been fixed for each $M_\alpha^{r_{\alpha}}$, the remaining exponents in the sequence
$\{k_{1\alpha}\} \cup \{k_{i\alpha}\}_{i=2}^\xi$ are constrained via
(\ref{eq:SumOfExponentsCondition}). We can recursively examine each
variable $x_i$ in turn using these constraints to build a set of
sequences
\begin{equation}
\label{eq:ProductSequenceDefinition}
S_l = \{S_{l\alpha}\}_{\alpha=1}^m = \{(k_{1\alpha},k_{2\alpha},\dots,k_{\xi\alpha})\}_{\alpha=1}^m
\end{equation} 

\noindent where each $S_{l\alpha}$ defines the exponent sequence for
its polynomial $M_{\alpha}^{r_{\alpha}}$ that
will produce the target term $T$ after the product is expanded. The
maximum value of $l$ depends on the target term $T$ and how many
possible $v_{i\alpha}$ values are filtered out using constraints
(\ref{eq:SumOfExponentsEqualsTarget}) and
(\ref{eq:SumOfExponentsCondition}) at each step in the recursion.

Once the set ${\bf S}=\{S_l\}$ has been constructed, we use
Equation (\ref{eq:MultinomialCoefficientEquation}) on each polynomial's
$\{k_{i\alpha}\}_{i=1}^\xi$ in $S_{l\alpha}$ to find the contributing coefficients. The
final coefficient value for term $T$ resulting from operation $G^i$ is
\begin{equation}
\label{eq:TargetTermValueForSingleGroupElement}
t_i = \sum_l \tau_l = \sum_l \prod_{\alpha=1}^m {d_\alpha \choose S_{l\alpha}}.
\end{equation}

\noindent To find the total number of unique colorings under the
group action, this process is applied to each element $G^i \in {\bf
  G}$ and the results are summed and then divided by $|{\bf G}|$.

We can further optimize the search for contributing terms by ordering
the exponents in the target term $T$ in descending order. Because the possible sequences
$\{k_{1\alpha}\}_{\alpha=1}^m$ are filtered using $T_1$, larger values
for $T_1$ are more likely to result in smaller sets of
$\{k_{i\alpha}\}_{\alpha=1}^m$ across the polynomials. All the
$\{k_{1\alpha}\}_{\alpha=1}^m$ need to sum to $T_1$
(\ref{eq:SumOfExponentsEqualsTarget}); if $T_1$ has 
smaller values (like 1 or 2), we will end up with lots of possible ways
to arrange them to sum to $T_1$ (which is not the the case for the larger
values). Since the final set of sequences $S_l$ is formed using a cartesian product, having a
few extra sequences from the $T_1$ pruning multiplies the total number
of sequences significantly. Additionally, constraint
(\ref{eq:SumOfExponentsCondition}) applied within each polynomial will
also reduce the total number of sequences to consider if the first
variables $x_1, x_2$, etc. are larger integers.

\subsection{Pseudocode Implementation}

Note. Implementations in python and Fortran are available in the supplementary material.

For both algorithms presented below, the operator ($\Leftarrow$) pushes
the value to its right onto the list to its left.

\begin{algorithm}
\caption{Recursive Sequence Constructor}\label{AlgSequenceConstructor}
\begin{algorithmic}
\Procedure{initialize}{$i$, $k_{i\alpha}$, $M_\alpha^{r_{\alpha}}$, $V_{\alpha}$,
  ${\bf T}$}
\State \emph{Constructs a Sequence Object tree recursively for a }
\State \emph{single $M_{\alpha}^{r_{\alpha}}$ by filtering possible exponents on each
  $x_i$}
\State \emph{in the polynomial. The object has the following }
\State \emph{properties:}
\State \; root: $k_{i\alpha}$, proposed exponent of variable $x_i$ in $M_\alpha^{r_{\alpha}}$.
\State \; parent: proposed Sequence object for $k_{i-1,\alpha}$ of $x_{i-1}$.
\State \; used: the sum of the proposed exponents to left of
\State \; \hspace{20pt}  and including this variable $\sum_{j=1}^i k_{i\alpha}$.
\State
\State $i$: index of variable in  $M_\alpha^{r_{\alpha}}$
\State $k_{i\alpha}$: proposed exponent of variable $x_i$ in $M_\alpha^{r_{\alpha}}$.
\State $M_\alpha^{r_{\alpha}}$: P\'olya polynomial representation of a
single 
\State \hspace{22pt} polynomial in $P^i$ (\ref{eq:ProductOfMultinomials}).
\State $V_{\alpha}$: possible exponents for $M_{\alpha}^{r_{\alpha}}$
(\ref{eq:PossibleVariableValues}).
\State ${\bf T}$: $\{T_i\}_{i=1}^{\xi}$ exponents for the
concentration term (\ref{eq:FixedColorTerm}). 
\State \dotfill
\If {$i = 1$}
  \State \textit{self}.used $\gets$ \textit{self}.root +
  \textit{self}.parent.used
\Else
  \State \textit{self}.used $\gets$ \textit{self}.root
\EndIf
\State \textit{self}.kids $\gets$ empty
\If {$i \le \xi$}
\For {$p \in V_\alpha$}
\State $rem \gets$ $p$ - \textit{self}.root
\If {$0 \le rem \le T_i$ {\bf and}
  $|rem| \le d_{\alpha}r_{\alpha}-\textit{self}.\textrm{used}$ \\ \pushcode[0]
  \pushcode[0] \pushcode[0] {\bf and}
  $|p-\textit{self}.\textrm{used}|\, \textrm{mod}\, r_{\alpha} = 0$}
\State $\textit{self}.\textrm{kids} \Leftarrow \textrm{Sequence}(i+1, rem,
M_{\alpha}^{r_{\alpha}}, V_{\alpha}, {\bf T})$
\EndIf
\EndFor
\EndIf
\EndProcedure \\ \hrulefill
\Function{expand}{sequence}
\State \emph{Generates a set of $S_{l\alpha}$ from a single Sequence
  object.}
\State sequence: the object created using \textsc{initialize}.
\State \dotfill
\State $\textit{sequences} \gets $ empty
\For {$\textit{kid} \in $ sequence.kids}
\For {$\textit{seq} \in \textsc{expand}(kid)$}
\State $\textit{sequences} \Leftarrow kid.\textrm{root} \cup \textit{seq}$
\EndFor
\EndFor
\If {$\textrm{len(sequence.kids)} = 0$}
  \State $\textit{sequences} \gets \{kid.\textrm{root}\}$
\EndIf
\State \Return {\textit{sequences}}
\EndFunction \\ \hrulefill
\Function{build\_S$_l$}{${\bf k}$, ${\bf V}$, $P^i$, ${\bf T}$}
\State \emph{Constructs $S_l$ from $\{k_{1\alpha}\}_{\alpha=1}^m$ for
  a $P^i$} (\ref{eq:ProductOfMultinomials}).
\State ${\bf k}$: $\{k_{1\alpha}\}_{\alpha=1}^m$ set of possible
exponent values on the 
\State \hspace{10pt} \emph{first} variable in each
$M_{\alpha}^{r_{\alpha}} \in P^i$.
\State ${\bf V}$: $\{V_{\alpha}\}_{\alpha=1}^m$ possible exponents for
each $M_{\alpha}^{r_{\alpha}}$ (\ref{eq:PossibleVariableValues}).
\State $P^i$: P\'olya polynomial representation
for a single
\State \hspace{14pt} operation in the group ${\bf G}$
(\ref{eq:ProductOfMultinomials}). 
\State ${\bf T}$: $\{T_i\}_{i=1}^{\xi}$ exponents for the
concentration term (\ref{eq:FixedColorTerm}).
\State \dotfill
\State $\textit{sequences} \gets \textrm{empty}$
\For {$\alpha \in \{1\dots m\}$}
  \State $\textit{seq} \gets \textsc{initialize}(1, k_{1\alpha},
  M_{\alpha}^{r_{\alpha}}, V_{\alpha}, {\bf T})$
  \State $\textit{sequences} \Leftarrow \textsc{expand}(seq)$
\EndFor
\State \Return $\textit{sequences}$
\EndFunction
\end{algorithmic}
\end{algorithm}

For algorithm (\ref{AlgSequenceConstructor}) in the \textsc{expand}
procedure, the $\cup$ operator horizontally concatenates
the integer \textit{root} to an existing sequence of integers.

For \textsc{build\_S$_l$}, we use the exponent $k_{1\alpha}$ on the first
variable in each polynomial to construct a full set of possible
sequences for that polynomial. Those sets of sequences are
then combined in \textsc{sum\_sequences}
(alg. \ref{AlgCoefficientCalculator}) using a cartesian product
over the sets in each multinomial.

\begin{algorithm}
\caption{Coefficient Calculator}\label{AlgCoefficientCalculator}
\begin{algorithmic}
\Function{sum\_sequences}{S$_l$}
\State \emph{Finds $\tau_{l}$} (\ref{eq:TargetTermValueForSingleGroupElement})
\emph{for $S_l= \{S_{l\alpha}\}_{\alpha=1}^m$} (\ref{eq:ProductSequenceDefinition})
\State S$_l$: a set of lists (of exponent sequences
$\{k_{i\alpha}\}_{i=1}^{\xi}$)
\State \hspace{10pt} for each polynomial $M_{\alpha}^{r_{\alpha}}$ in the product $P^i$
(\ref{eq:ProductOfMultinomials}).
\State \dotfill
\State $K_l \gets S_{l1}\times
S_{l2} \times \dots \times S_{lm} = \langle
\{(k_{i\alpha})_{i=1}^\xi\}_{\alpha=1}^m \rangle_l$
\State $\textit{coeff} \gets 0$
\For {{\bf each} $\{(k_{i\alpha})_{i=1}^\xi\}_{\alpha=1}^m \in K_l$}
\If {$\sum_{\alpha=1}^m k_{i\alpha} = T_i \; \forall \; i \in \{1\dots \xi\}$}
\State $\textit{coeff} \gets \textit{coeff} + \prod_{\alpha=1}^m {d_\alpha \choose \{k_{i\alpha}\}_{i=1}^{\xi}}$
\EndIf
\EndFor
\State \Return $\textit{coeff}$
\EndFunction \\ \hrulefill
\Function{coefficient}{${\bf T}$, $P^i$, ${\bf V}$}
\State \emph{Constructs } ${\bf S} = \{S_l\}$ \emph{and
  calculates $t_i$} (\ref{eq:TargetTermValueForSingleGroupElement})
\State ${\bf T}$: $\{T_i\}_{i=1}^{\xi}$ exponents for the concentration term (\ref{eq:FixedColorTerm}).
\State $P^i$: P\'olya polynomial representation
for a single
\State \hspace{14pt} operation in the group ${\bf G}$
(\ref{eq:ProductOfMultinomials}). 
\State ${\bf V}$: $\{V_{\alpha}\}_{\alpha=1}^m$ possible exponents for
each $M_{\alpha}^{r_{\alpha}}$ (\ref{eq:PossibleVariableValues}).
\State \dotfill
\If {$m=1$}
  \If {$r_1>T_i \; \forall \; i=1..\xi$}
    \State \Return $0$
  \Else
    \State \Return ${d_1 \choose {T_1T_2\dots T_\xi}}$
  \EndIf
\Else
  \State ${\bf T} \gets \textrm{sorted}({\bf T})$
  \State $\textit{possible} \gets V_1 \times V_2 \times \dots \times
  V_m$
  \State $\textit{coeffs} \gets 0$
  \For {$\{k_{1\alpha}\}_{\alpha=1}^m \in \textit{possible}$}
    \If {$\sum_{\alpha=1}^m k_{1\alpha} = T_1$}
      \State $S_l \gets
      \textsc{build\_S}_l(\{k_{1\alpha}\}_{\alpha=1}^m, {\bf V}, P^i, {\bf T})$
      \State $\textit{coeffs} \gets \textit{coeffs} + \textsc{sum\_sequences}(S_l)$
    \EndIf
  \EndFor
  \State \Return \textit{coeffs}
\EndIf
\EndFunction
\end{algorithmic}
\end{algorithm}

For algorithm (\ref{AlgCoefficientCalculator}) in the
\textsc{sum\_sequences} function, $K_l$ is calculated
using the cartesian product of the individual $S_{l\alpha}$, where for a
given $l$, the number of sequences $\{k_i\}_{i=1}^\xi \in S_{l\alpha}$ may be
arbitrary. For example, a product of three polynomials
$M_1^4M_2^3M_3^2$ may produce possible sequences with $|S_{l1}|=2$,
$|S_{l2}|=4$ and $|S_{l3}|=4$. Then $|K_l|=2 \times 4 \times 4=32$ and
each element in $K_l$ is a set of three sequences:
$\{(k_{i\alpha})_{i=1}^\xi\}_{\alpha=1}^3$, one for each polynomial, which
specifies the exponents on the contributing term from that
polynomial. Also, when calculating multinomial coefficients, we use
the form in eqn. (\ref{eq:MultinomialCoefficientEquation}) in terms of
binomial coefficients with a fast, stable algorithm from Manolopoulos
\cite{Manolo2002}.

In practice, many of the group operations $G^i$ produce identical
products $M_1^{r_1}M_2^{r_2}\dots M_m^{r_m}$. Thus before computing
any of the coefficients from the polynomials, we first form the
polynomial products for each group operation and then add identical
products together.

\section{Computational Order and Performance}

The algorithm is structured around the \emph{a priori} knowledge of
the fixed concentration term (\ref{eq:FixedColorTerm}). At the
earliest possibility, we prune terms from individual polynomials that
would not contribute to the final polya coefficient in the expanded
product of polynomials. Because the P\'olya polynomial for each group
operation is based on its disjoint-cyclic form, the complexity of the
search can vary drastically from one group operation to the next. That
said, it is common for groups to have several classes whose group
operations (within each class) will have similar disjoint-cyclic forms
and thus also scale similarly. However, from group to group, the set
of classes and disjoint-cyclic forms may be very different; this makes
it difficult to make a statement about the scaling of the algorithm in
general. Although we could make statements about the scaling of
well-known sets of groups (for example the dihedral groups used in our
example above), we decided instead to craft certain special groups
with specific properties and run tests to determine the scaling
numerically.

\begin{figure}[]
\centerline{\includegraphics[width=.5\textwidth]{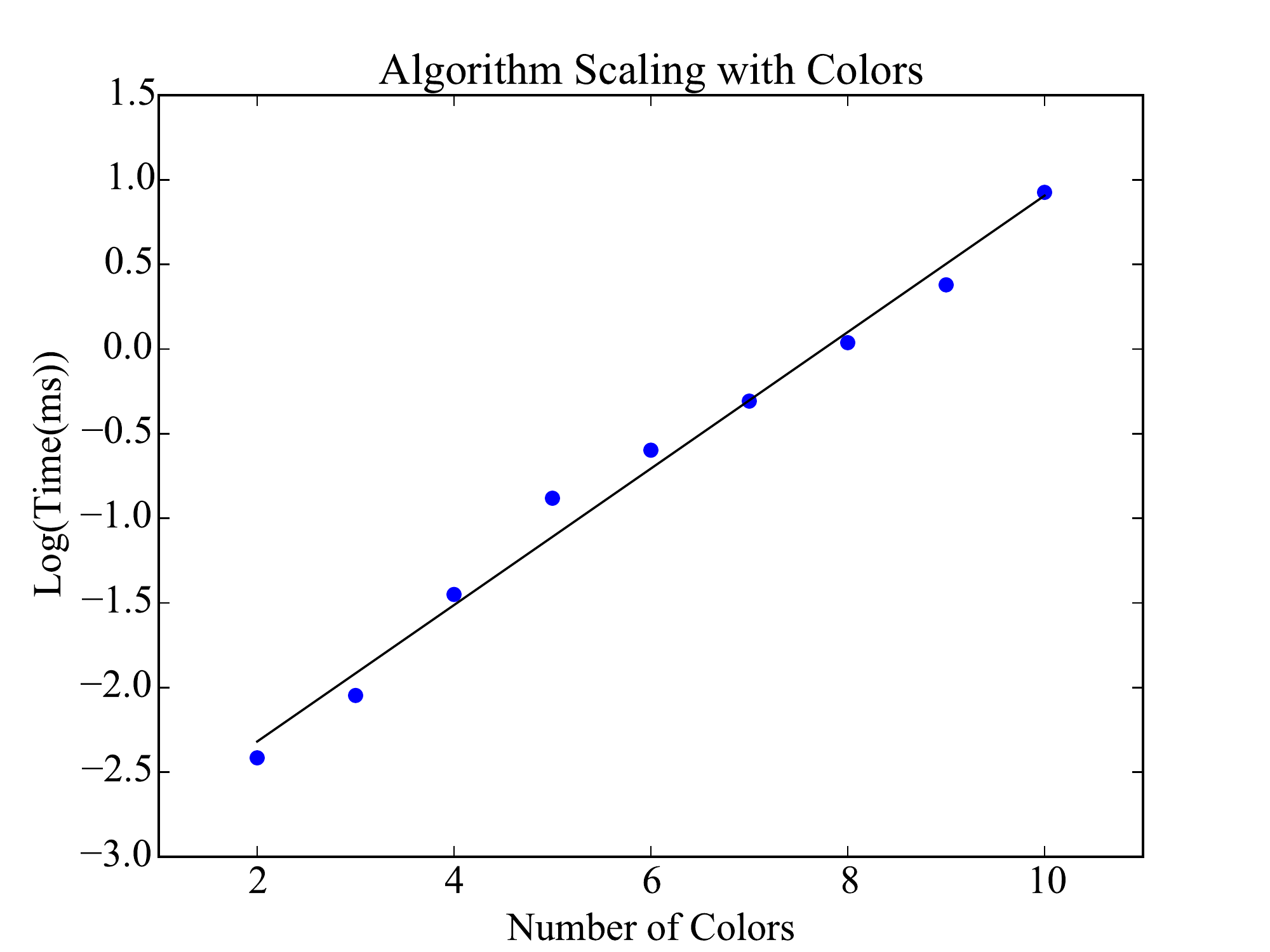}}
\caption[]{\label{fig:AlgorithmScalingColors} Log plot of the
  algorithm scaling as the number of colors increases. Since the
  number of variables $x_i$ in each polynomial increases with the
  number of colors, the combinatoric complexity of the expanded
  polynomial increases drastically with each additional color; this
  leads to an exponential scaling. The linear fit to the logarithmic
  data has a slope of 0.403.}
\end{figure}

In Figure \ref{fig:AlgorithmScalingColors} we plot the algorithm's
scaling as the number of colors in the enumeration increases. For each
$r$-cycle in the disjoint-cyclic form of a group operation, we
construct a polynomial with $\xi$ variables, where $\xi$ is the number of
colors used in the enumeration. Because the group operation results in
a product of these polynomials, increasing the number of colors by 1
increases the combinatoric complexity of the polynomial
\emph{expansion} exponentially. For this scaling experiment, we used
the same transitive group acting on a finite set with 20 elements for
each data point, but increased the number of colors in the fixed color
term $T$. We chose $T$ by dividing the number of elements in the group
as equally as possible; thus for 2 colors, we used $[10,10]$; for 3
colors we used $[8,6,6]$, then $[5,5,5,5]$, $[4,4,4,4,4]$, etc. Figure
\ref{fig:AlgorithmScalingColors} plots the $\log_{10}$ of the
execution time (in ms) as the number of colors increases. As expected,
the scaling is linear (on the log plot). The linear fit to the data
points has a slope of $0.403$.

\begin{figure}[htbp]
\centerline{\includegraphics[width=.5\textwidth]{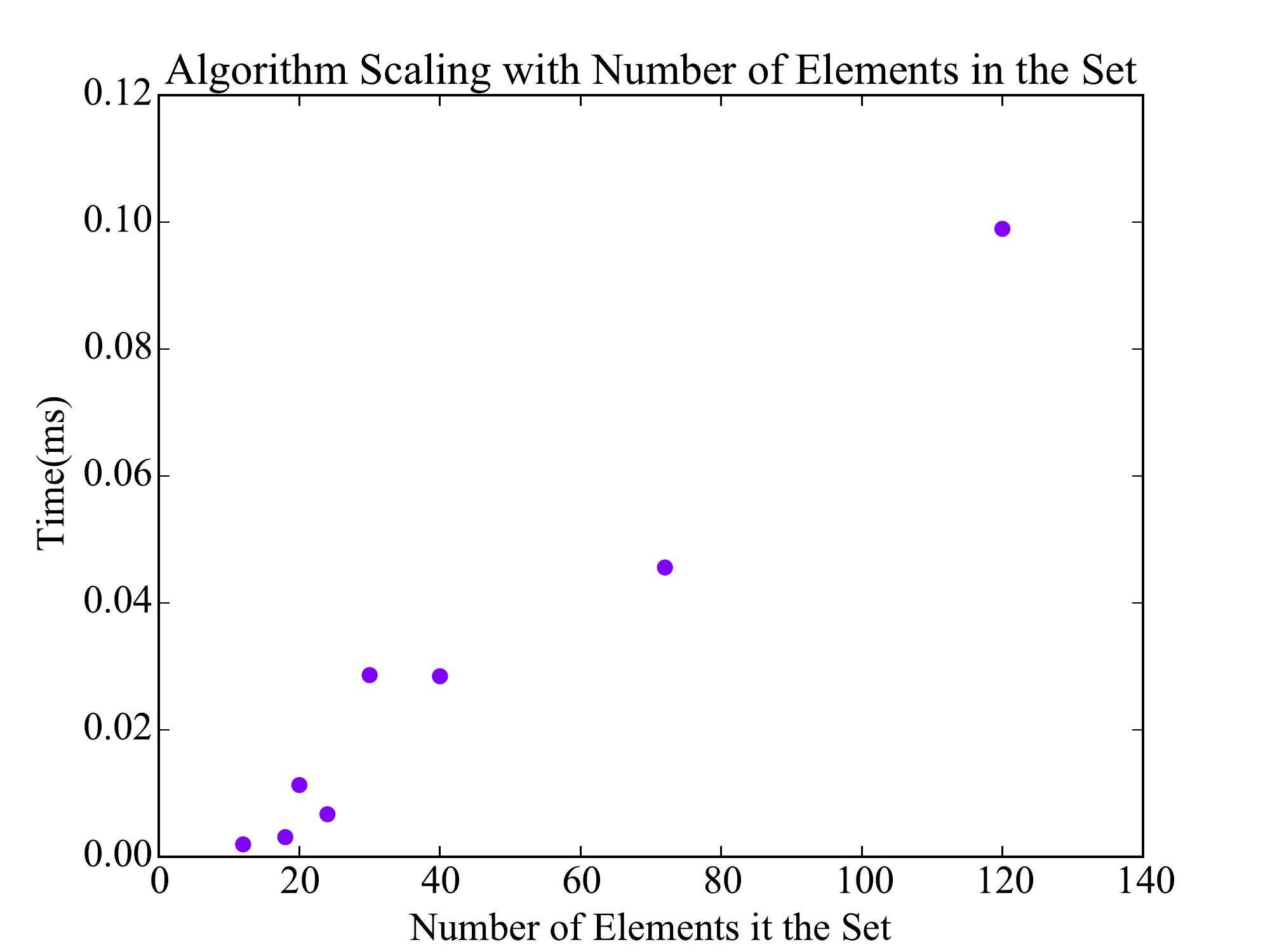}}
\caption[]{\label{fig:AlgorithmScalingSites} Algorithm scaling as the
  number of elements in the finite set increases. The P\'olya
  polynomial arises from the group operations' disjoint-cyclic form,
  so that more elements in the set results in a richer spectrum of
  possible polynomials multiplied together. Because of the algorithms
  aggresive pruning of terms, the exact disjoint-cyclic form of
  individual group operations has a large bearing on the algorithm's
  scaling. As such it isn't surprising that there is some scatter in
  the timings as the number of elements in the set increases.}
\end{figure}

As the number of elements in the finite set increases, the possible P\'olya
polynomial representations for each group operation's disjoint-cyclic
form increases exponentially. In the worst case, a group acting on a
set with $k$ elements may have an operation with $k$ 1-cycles; on the
other hand, that same group may have an operation with a single
$k$-cycle, with lots of possibilities in between. Because of the
richness of possibilities, it is almost impossible to make general
statements about the algorithm's scaling without knowing the structure
of the group and its classes. In Figure
\ref{fig:AlgorithmScalingSites}, we plot the scaling for a set of
related groups (all are isomorphic to the direct product of S$_3$ $\times$ S$_4$)
applied to finite sets of varying sizes. Every data point was 
generated using a transitive group with 144 elements. Thus, this plot
shows the algorithm's scaling when the group is the same and the
number of elements in the finite set changes. Although the scaling
appears almost linear, there is a lot of scatter in the data. Given
the rich spectrum of possible P\'olya polynomials that we can form as
the set size increases, the scatter isn't surprising.

\begin{figure}[htbp]
\centerline{\includegraphics[width=.5\textwidth]{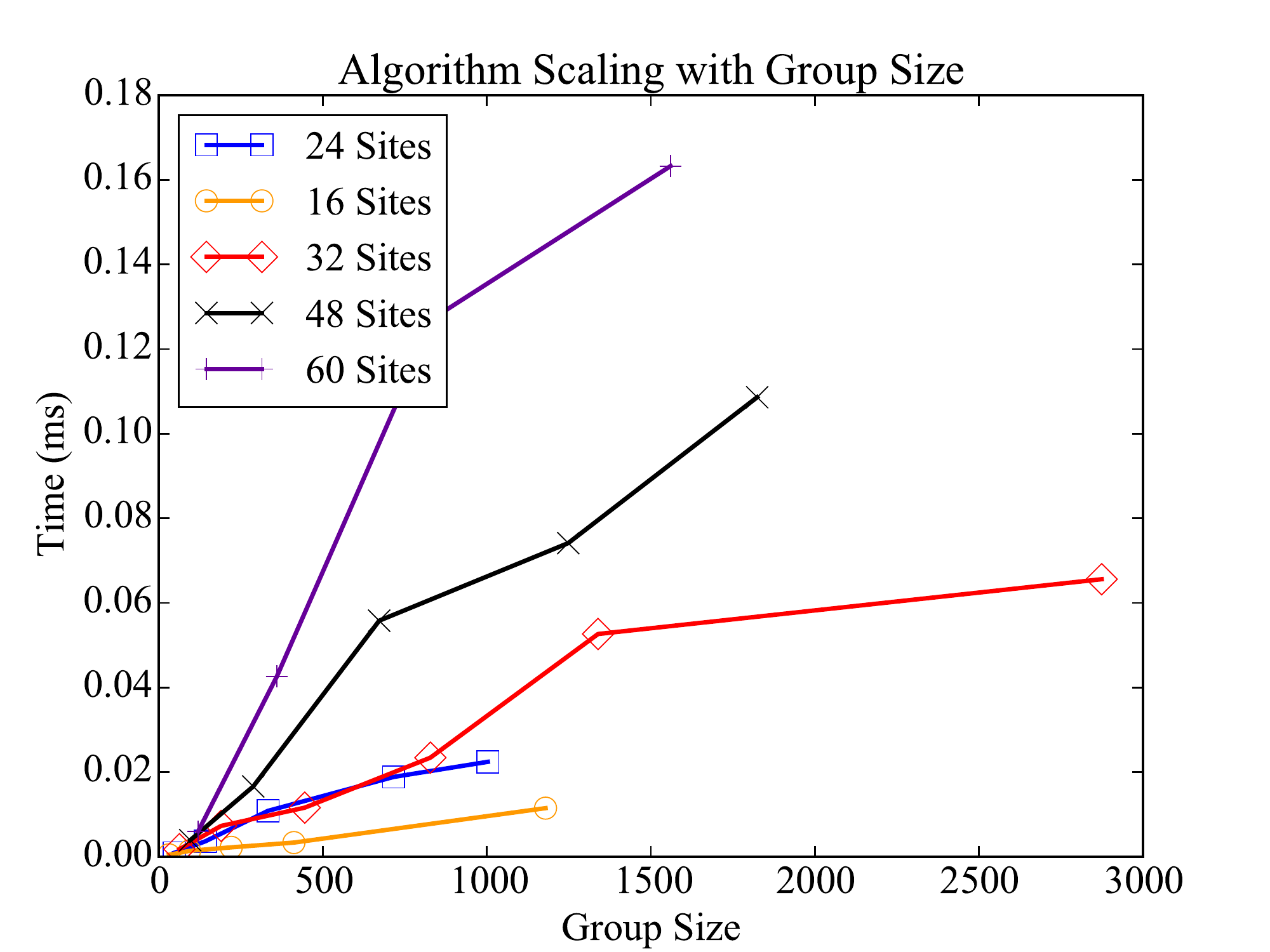}}
\caption[]{\label{fig:AlgorithmScalingGroupSize} Algorithm scaling
  with group size for an enumeration problem from solid state physics
  \cite{Hart:2008}. We used the unique permutation groups arising from
  all derivative super structures of a simple cubic lattice for a
  given number of sites in the unit cell. The behavior is generally
  linear with increasing group size. }
\end{figure}

Finally, we consider the scaling as the group size increases. For this
test, we selected the set of unique groups arising from the
enumeration of all derivative super structures of a simple cubic
lattice for a given number of sites in the unit cell
\cite{Hart:2008}. Since the groups are formed from the symmetries of
real crystals, they arise from the semidirect product of operations
related to physical rotations and translations of the crystal. In this
respect, they have similar structure for comparison. In most cases,
the scaling is obviously linear; however, the slope of each trend
varies from group to group. This once again
highlights the scaling's heavy dependence on the specific
disjoint-cyclic forms of the group operations. Even for groups with
obvious similarity, the scaling may be different.

\section{Summary}

Until now, no low-level, numerical implementation of P\'olya's
enumeration theorem was readily available; instead, a computer algebra
system (CAS) was used to symbolically solve the polynomial expansion
problem posed by P\'olya. While such systems are effective for small,
simpler calculations, as the difficulty of the problem increases, they
become impractical solutions. Additionally, codes that perform the
actual enumeration of the colorings are often implemented in low-level
codes and interoperability with a CAS is not necessarily easy to automate.

We presented a low-level, purely numerical algorithm that exploits the
properties of polynomials to restrict the combinatoric complexity of
the expansion. By considering only those coefficients in the unexpanded
polynomials that might contribute to the final answer, the algorithm
reduces the number of terms that must be included to find the
significant term in the expansion. 

Because of the algorithm scaling's reliance on the exact structure of
the group and the disjoint-cyclic form of its operations, a rigorous
analysis of the scaling is not possible without knowledge of the
group. Instead, we presented some numerical timing results from
representative, real-life problems that show the general scaling
behavior. Because all the timings are in the millisecond to second regime
anyway, a more rigorous analysis of the algorithm's scaling is unnecessary.

In contrast to the CAS solutions whose execution times range from
milliseconds to hours, our algorithm consistently performs in the
millisecond to second regime, even for complex problems. Additionally,
it is easy to implement in low-level languages, making it useful for
confirming enumeration results. This makes it an effective
substitute for alternative CAS implementations.

In computational materials science, chemistry, and related subfields such as
  computational drug discovery, combinatorial searches are becoming
  increasingly important, especially in high-throughput
  studies \cite{nmat_review}. The upside potential of these efforts
  continues to grow because computing power continues to become cheaper
  and algorithms continue to evolve. As computational methods gain a
  larger market share in materials discovery, algorithms such as this
  one are important as they provide validation support to complex
  simulation codes. The present algorithm has been useful in checking
  a new algorithm extending the work in
  Refs.~\cite{Hart:2008,Hart:2009,Hart:2012}, and P\'olya's theorem was
  recently used in Mustapha's enumeration 
  algorithm\cite{Mustapha:2013} that has been incorporated into the CRYSTAL14 software
  package \cite{QUA:QUA24658}.
\begin{acknowledgements}
This work was supported under ONR (MURI N00014-13-1-0635). 

\end{acknowledgements}

\bibliographystyle{spmpsci}      

\section{Supplementary Material}

The source code to implement this algorithm is available for both
python and Fortran at:\\

https://github.com/rosenbrockc/polya \\

\noindent The home page on github has full
instructions for using either version of the code as well a battery of
over 50 unit tests that were used to verify and time the
algorithm. The unit tests can be executed using the \textsc{fortpy}
framework available via the Python Package Index. Instructions for
running the unit tests are also on the github home page.

\end{document}